\documentclass[11pt]{article}
\usepackage{amssymb}
\usepackage[all]{xy}
\usepackage[latin1]{inputenc}
\usepackage{amsfonts}
\usepackage{amsmath}
\usepackage{amssymb}
\usepackage{url}
\usepackage{hyperref}
\usepackage{graphicx}

\providecommand{\U}[1]{\protect\rule{.1in}{.1in}}
\newtheorem{theo}{Theorem}[section]
\newtheorem{prop}[theo]{Proposition}
\newtheorem{lem}[theo]{Lemma}
\newtheorem{cor}[theo]{Corollary}

\newtheorem{defi}[theo]{Definition}

\numberwithin{equation}{section}

\begin{document}

\title{A Backward Particle Interpretation of Feynman-Kac Formulae}
\author{ Pierre Del Moral\thanks{%
Centre INRIA Bordeaux et Sud-Ouest \& Institut de Mathématiques de Bordeaux
, Université de Bordeaux I, 351 cours de la Libération 33405 Talence cedex,
France (Pierre.Del-Moral@inria.fr)} , Arnaud Doucet\thanks{%
Department of Statistics \& Department of Computer Science, University of
British Columbia, 333-6356 Agricultural Road, Vancouver, BC, V6T 1Z2, Canada
and The Institute of Statistical Mathematics, 4-6-7 Minami-Azabu, Minato-ku,
Tokyo 106-8569, Japan (Arnaud@stat.ubc.ca)} , Sumeetpal S. Singh\thanks{%
Department of Engineering, University of Cambridge, Trumpington Street, CB2
1PZ, United Kingdom (sss40@cam.ac.uk)}}
\date{20th July 2009 \\
(Preprint submitted to ESAIM: M2AN)}
\maketitle

\begin{abstract}
We design a particle interpretation of Feynman-Kac measures on path spaces
based on a backward Markovian representation combined with a traditional
mean field particle interpretation of the flow of their final time
marginals. In contrast to traditional genealogical tree based models, these
new particle algorithms can be used to compute normalized additive
functionals \textquotedblleft on-the-fly\textquotedblright\ as well as their
limiting occupation measures with a given precision degree that does not
depend on the final time horizon.

We provide uniform convergence results w.r.t. the time horizon parameter as
well as functional central limit theorems and exponential concentration
estimates. We also illustrate these results in the context of computational
physics and imaginary time Schroedinger type partial differential equations,
with a special interest in the numerical approximation of the invariant
measure associated to $h$-processes.\newline


\end{abstract}

\section{Introduction}

Let $(E_{n})_{n\geq 0}$ be a sequence of measurable spaces equipped with
some $\sigma $-fields $(\mathcal{E}_{n})_{n\geq 0}$, and we let $\mathcal{P}%
(E_{n})$ be the set of all probability measures over the set $E_{n}$, with $%
n\geq 0$. We let $X_{n}$ be a Markov chain with Markov transition $M_{n}$ on
$E_{n}$, and we consider a sequence of $(0,1]$-valued potential functions $%
G_{n}$ on the set $E_{n}$. The Feynman-Kac path measure associated with the
pairs $(M_{n},G_{n})$ is the probability measure $\mathbb{Q}_{n}$ on the
product state space $E_{[0,n]}:=\left( E_{0}\times \ldots \times
E_{n}\right) $ defined by the following formula
\begin{equation}
d\mathbb{Q}_{n}:=\frac{1}{\mathcal{Z}_{n}}~\left\{ \prod_{0\leq
p<n}G_{p}(X_{p})\right\} ~d\mathbb{P}_{n}  \label{defQn}
\end{equation}%
where $\mathcal{Z}_{n}$ is a normalizing constant and $\mathbb{P}_{n}$ is
the distribution of the random paths $(X_{p})_{0\leq p\leq n}$ of the Markov
process $X_{p}$ from the origin $p=0$, up to the current time $p=n$. We also
denote by $\Gamma _{n}=\mathcal{Z}_{n}~\mathbb{Q}_{n}$ its unnormalized
version.

These distributions arise in a variety of application areas, including
filtering, Bayesian inference, branching processes in biology, particle
absorption problems in physics and many other instances. We refer the reader
to the pair of books~\cite{fk,doucet} and references therein. Feynman-Kac
models also play a central role in the numerical analysis of certain partial
differential equations, offering a natural way to solve these functional
integral models by simulating random paths of stochastic processes. These
Feynman-Kac models were originally presented by Mark Kac in 1949~\cite{kac}
for continuous time processes. These continuous time models are used in
molecular chemistry and computational physics to calculate the ground state
energy of some Hamiltonian operators associated with some potential function
$V$ describing the energy of a molecular configuration (see for instance~%
\cite{benj2,dms,benjamin,rousset}, and references therein).

To better connect these partial differential equation models with (\ref%
{defQn}), let us assume that $M_{n}(x_{n-1},dx_{n})$ is the Markov
probability transition $X_{n}=x_{n}\leadsto X_{n+1}=x_{n+1}$ coming from a
discretization in time $X_{n}=X_{t_{n}}^{\prime }$ of a continuous time $E$%
-valued Markov process $X_{t}^{\prime }$ on a given time mesh $%
(t_{n})_{n\geq 0}$ with a given time step $(t_{n}-t_{n-1})=\Delta t$. For
potential functions of the form $G_{n}=e^{-V\Delta t}$, the measures $%
\mathbb{Q}_{n}\simeq _{\Delta t\rightarrow 0}\mathbb{Q}_{t_{n}}$ represents
the time discretization of the following distribution:
\begin{equation*}
d\mathbb{Q}_{t}=\frac{1}{\mathcal{Z}_{t}}~\exp {\left(
-\int_{0}^{t}~V(X_{s}^{\prime })~ds\right) }~d\mathbb{P}_{t}^{X^{\prime }}
\end{equation*}%
where $\mathbb{P}_{t}^{X^{\prime }}$ stands for the distribution of the
random paths $(X_{s}^{\prime })_{0\leq s\leq t}$ with a given infinitesimal
generator $L$. The marginal distributions $\gamma _{t}$ at time $t$ of the
unnormalized measures $\mathcal{Z}_{t}~d\mathbb{Q}_{t}$ are the solution of
the so-called imaginary time Schroedinger equation, given in weak
formulation on every sufficiently regular function $f$ by
\begin{equation*}
\frac{d}{dt}~\gamma _{t}(f):=\gamma _{t}(L^{V}(f))\quad \mbox{\rm with}\quad
L^{V}=L-V
\end{equation*}%
The errors introduced by the discretization of the time are well understood
for regular models, we refer the interested reader to~\cite%
{delmojacod,dimasi,korez2,picard} in the context of nonlinear filtering.

In this article, we design an numerical approximation of the distributions $%
\mathbb{Q}_{n}$ based on the simulation of a sequence of mean field
interacting particle systems. In molecular chemistry, these evolutionary
type models are often interpreted as a quantum or diffusion Monte Carlo
model. In this context, particles often are referred as walkers, to
distinguish the virtual particle-like objects to physical particles, like
electrons of atoms. In contrast to traditional genealogical tree based
approximations (see for instance~\cite{fk}), the particle model presented in
this article can approximate additive functionals of the form
\begin{equation}
\overline{F}_{n}(x_{0},\ldots ,x_{n})=\frac{1}{(n+1)}\sum_{0\leq p\leq
n}f_{p}(x_{p})  \label{additivefunctnorm}
\end{equation}%
uniformly with respect to the time horizon. Moreover this computation can be
done \textquotedblleft on-the-fly\textquotedblright . To give a flavor of
the impact of these results, we recall that the precision of the algorithm
corresponds to the size $N$ of the particle system. If $\mathbb{Q}_{n}^{N}$
stands for the $N$-particle approximation of $\mathbb{Q}_{n}$, under some
appropriate regularity properties, we shall prove the following uniform and
non asymptotic Gaussian concentration estimates\footnote{%
Consult the last paragraph of this section for a statement of the notation
used in this article.}:
\begin{equation*}
{\frac{1}{N}\log {\sup_{n\geq 0}{\mathbb{P}\left( \left\vert [\mathbb{Q}%
_{n}^{N}-\mathbb{Q}_{n}](F_{n})\right\vert \geq \frac{b}{\sqrt{N}}+\epsilon
\right) }}}\leq -{\epsilon ^{2}}/{(2b^{2})}
\end{equation*}%
for any $\epsilon >0$, and for some finite constant $b<\infty $. In the
filtering context, $\mathbb{Q}_{n}^{N}$ corresponds to the sequential\ Monte
Carlo approximation of the forward filtering backward smoothing recursion.
Recently, a theoretical study of this problem was undertaken by \cite%
{Douc2009}. Our results complement theirs and we
present functional central limit theorems as well as non-asymptotic variance bounds.
Additionally, we show how the forward filtering backward smoothing estimates of additive functionals can be computed using a forward only recursion. This has applications to online parameter estimation for non-linear non-Gaussian state-space models.

For time homogeneous models $(M_{n},f_{n},G_{n})=(M,f,G)$ associated with a
lower bounded potential function $G>\delta $, and a $M$-reversible
transition w.r.t. to some probability measure $\mu $ s.t. $M(x,\mbox{\LARGE
.})\sim \mu $ and $(M(x,\mbox{\LARGE .})/d\mu )\in \mathbb{L}_{2}(\mu )$, it
can be established that $\mathbb{Q}_{n}(F_{n})$ converges to $\mu _{h}(f)$,
as $n\rightarrow \infty $, with the measure $\mu _{h}$ defined below
\begin{equation*}
\mu _{h}(dx):=\frac{1}{\mu (hM(h))}~h(x)~M(h)(x)~\mu (dx)
\end{equation*}%
In the above display, $h$ is a positive eigenmeasure associated with the top
eigenvalue of the integral operator $Q(x,dy)=G(x)M(x,dy)$ on $\mathbb{L}%
_{2}(\mu )$ (see for instance section 12.4 in~\cite{fk}). This measure $\mu
_{h}$ is in fact the invariant measure of the $h$-process defined as the
Markov chain $X^{h}$ with elementary Markov transitions $M_{h}(x,dy)\propto
M(x,dy)h(y)$. As the initiated reader would have certainly noticed, the
above convergence result is only valid under some appropriate mixing
conditions on the $h$-process. The long time behavior of these $h$-processes
and their connections to various applications areas of probability,
analysis, geometry and partial differential equations, have been the subject
of countless papers for many years in applied probability. In our framework,
using elementary manipulations, the Gaussian estimate given above can be
used to calibrate the convergence of the particle estimate $\mathbb{Q}%
_{n}^{N}(F_{n})$ towards $\mu _{h}(f)$, as the pair of parameters $N$ and $%
n\rightarrow \infty $.

The rest of this article is organized as follows:

In section~\ref{descript}, we describe the mean field particle models used
to design the particle approximation measures $\mathbb{Q}_{n}^{N}$. In
section~\ref{state}, we state the main results presented in this article,
including a functional central limit theorem, and non asymptotic mean error
bounds. Section~\ref{bacw} is dedicated to a key backward Markov chain
representation of the measures $\mathbb{Q}_{n}$. The analysis of our
particle approximations is provided in section~\ref{secpartic}. The final
two sections, section~\ref{fluctsec} and section~\ref{nonasymp}, are mainly
concerned with the proof of the two main theorems presented in section~\ref%
{state}.

For the convenience of the reader, we end this introduction with some
notation used in the present article. We denote respectively by $\mathcal{M}%
(E)$, and $\mathcal{B}(E)$, the set of all finite signed measures on some
measurable space $(E,\mathcal{E})$, and the Banach space of all bounded and
measurable functions $f$ equipped with the uniform norm $\Vert f\Vert $. We
let $\mu (f)=\int ~\mu (dx)~f(x)$, be the Lebesgue integral of a function $%
f\in \mathcal{B}(E)$, with respect to a measure $\mu \in \mathcal{M}(E)$. We
recall that a bounded integral kernel $M(x,dy)$ from a measurable space $(E,%
\mathcal{E})$ into an auxiliary measurable space $(E^{\prime },\mathcal{E}%
^{\prime })$ is an operator $f\mapsto M(f)$ from $\mathcal{B}(E^{\prime })$
into $\mathcal{B}(E)$ such that the functions
\begin{equation*}
x\mapsto M(f)(x):=\int_{E^{\prime }}M(x,dy)f(y)
\end{equation*}%
are $\mathcal{E}$-measurable and bounded, for any $f\in \mathcal{B}%
(E^{\prime })$. In the above displayed formulae, $dy$ stands for an
infinitesimal neighborhood of a point $y$ in $E^{\prime }$. The kernel $M$
also generates a dual operator $\mu \mapsto \mu M$ from $\mathcal{M}(E)$
into $\mathcal{M}(E^{\prime })$ defined by $(\mu M)(f):=\mu (M(f))$. A
Markov kernel is a positive and bounded integral operator $M$ with $M(1)=1$.
Given a pair of bounded integral operators $(M_{1},M_{2})$, we let $%
(M_{1}M_{2})$ the composition operator defined by $%
(M_{1}M_{2})(f)=M_{1}(M_{2}(f))$. For time homogenous state spaces, we
denote by $M^{m}=M^{m-1}M=MM^{m-1}$ the $m$-th composition of a given
bounded integral operator $M$, with $m\geq 1$. Given a positive function $G$
on $E$, we let $\Psi _{G}~:~\eta \in \mathcal{P}(E)\mapsto \Psi _{G}(\eta
)\in \mathcal{P}(E)$, be the Boltzmann-Gibbs transformation defined by
\begin{equation*}
\Psi _{G}(\eta )(dx):=\frac{1}{\eta (G)}~G(x)~\eta (dx)
\end{equation*}

\section{Description of the models}

\label{descript} The numerical approximation of the path-space distributions
(\ref{defQn}) requires extensive calculations. The mean field particle
interpretation of these models are based on the fact that the flow of the $n$%
-th time marginals $\eta _{n}$ of the measures $\mathbb{Q}_{n}$ satisfy a
non linear evolution equation of the following form
\begin{equation}
\eta _{n+1}(dy)=\int \eta _{n}(dx)K_{n+1,\eta _{n}}(x,dy)  \label{nonlin1}
\end{equation}%
for some collection of Markov transitions $K_{n+1,\eta }$, indexed by the
time parameter $n\geq 0$ and the set of probability measures $\mathcal{P}%
(E_{n})$. The mean field particle interpretation of the nonlinear measure
valued model (\ref{nonlin1}) is the $E_{n}^{N}$-valued Markov chain
\begin{equation*}
\xi _{n}=\left( \xi _{n}^{1},\xi _{n}^{2},\ldots ,\xi _{n}^{N}\right) \in
E_{n}^{N}
\end{equation*}%
with elementary transitions defined as
\begin{equation}
\mathbb{P}\left( \xi _{n+1}\in dx~|~\xi _{n}\right)
=\prod_{i=1}^{N}~K_{n+1,\eta _{n}^{N}}(\xi _{n}^{i},dx^{i})\quad
\mbox{\rm
with}\quad \eta _{n}^{N}:=\frac{1}{N}\sum_{j=1}^{N}~\delta _{\xi _{n}^{j}}
\label{meanfield}
\end{equation}%
In the above displayed formula, $dx$ stands for an infinitesimal
neighborhood of the point $x=(x^{1},\ldots ,x^{N})\in E_{n+1}^{N}$. The
initial system $\xi _{0}$ consists of $N$ independent and identically
distributed random variables with common law $\eta _{0}$. We let $\mathcal{F}%
_{n}^{N}:=\sigma \left( \xi _{0},\ldots ,\xi _{n}\right) $ be the natural
filtration associated with the $N$-particle approximation model defined
above. The resulting particle model coincides with a genetic type stochastic
algorithm $\xi _{n}\leadsto \widehat{\xi }_{n}\leadsto \xi _{n+1}$ with
selection transitions $\xi _{n}\leadsto \widehat{\xi }_{n}$ and mutation
transitions $\widehat{\xi }_{n}\leadsto \xi _{n+1}$ dictated by the
potential (or fitness) functions $G_{n}$ and the Markov transitions $M_{n+1}$%
.

During the selection stage $\xi _{n}\leadsto \widehat{\xi }_{n}$, for every
index $i$, with a probability $\epsilon _{n}G_{n}(\xi _{n}^{i})$, we set $%
\widehat{\xi }_{n}^{i}=\xi _{n}^{i}$, otherwise we replace $\xi _{n}^{i}$
with a new individual $\widehat{\xi }_{n}^{i}=\xi _{n}^{j}$ randomly chosen
from the whole population with a probability proportional to $G_{n}(\xi
_{n}^{j})$. The parameter $\epsilon _{n}\geq 0$ is a tuning parameter that
must satisfy the constraint $\epsilon _{n}G_{n}(\xi _{n}^{i})\leq 1$, for
every $1\leq i\leq N$. For $\epsilon _{n}=0$, the resulting proportional
selection transition corresponds to the so-called simple genetic model.
During the mutation stage, the selected particles $\widehat{\xi }%
_{n}^{i}\leadsto \xi _{n+1}^{i}$ evolve independently according to the
Markov transitions $M_{n+1}$.

If we interpret the selection transition as a birth and death process, then
arises the important notion of the ancestral line of a current individual.
More precisely, when a particle $\widehat{\xi }_{n-1}^{i}\longrightarrow \xi
_{n}^{i}$ evolves to a new location $\xi _{n}^{i}$, we can interpret $%
\widehat{\xi }_{n-1}^{i}$ as the parent of $\xi _{n}^{i}$. Looking backwards
in time and recalling that the particle $\widehat{\xi }_{n-1}^{i}$ has
selected a site $\xi _{n-1}^{j}$ in the configuration at time $(n-1)$, we
can interpret this site $\xi _{n-1}^{j}$ as the parent of $\widehat{\xi }%
_{n-1}^{i}$ and therefore as the ancestor denoted $\xi _{n-1,n}^{i}$ at
level $(n-1)$ of $\xi _{n}^{i}$. Running backwards in time we may trace the
whole ancestral line
\begin{equation}
\xi _{0,n}^{i}\longleftarrow \xi _{1,n}^{i}\longleftarrow \ldots
\longleftarrow \xi _{n-1,n}^{i}\longleftarrow \xi _{n,n}^{i}=\xi _{n}^{i}
\label{ancestraline}
\end{equation}%
More interestingly, the occupation measure of the corresponding $N$%
-genealogical tree model converges as $N\rightarrow \infty $ to the
conditional distribution $\mathbb{Q}_{n}$. For any function $F_{n}$ on the
path space $E_{[0,n]}$, we have the following convergence (to be stated
precisely later) as $N\rightarrow \infty $,
\begin{equation}
\lim_{N\rightarrow \infty }\frac{1}{N}\sum_{i=1}^{N}F_{n}(\xi _{0,n}^{i},\xi
_{1,n}^{i},\ldots ,\xi _{n,n}^{i})=\int ~\mathbb{Q}_{n}(d(x_{0},\ldots
,x_{n}))~F_{n}(x_{0},\ldots ,x_{n})  \label{cvarbre}
\end{equation}%
This convergence result can be refined in various directions. Nevertheless,
the asymptotic variance $\sigma _{n}^{2}(F_{n})$ of the above occupation
measure around $\mathbb{Q}_{n}$ increases quadratically with the final time
horizon $n$ for additive functions of the form
\begin{equation}
F_{n}(x_{0},\ldots ,x_{n})=\sum_{0\leq p\leq n}f_{p}(x_{p})\Rightarrow
\sigma _{n}^{2}(F_{n})\simeq n^{2}  \label{additivefunct}
\end{equation}%
with some collection of non negative functions $f_{p}$ on $E_{p}$. To be
more precise, let us examine a time homogeneous model $%
(E_{n},f_{n},G_{n},M_{n})=(E,f,G,M)$ with constant potential functions $%
G_{n}=1$ and mutation transitions $M$ s.t. $\eta _{0}M=\eta _{0}$. For the
choice of the tuning parameter $\epsilon =0$, using the asymptotic variance
formulae in~\cite[eqn. (9.13), page 304 ]{fk}, for any function $f$ s.t. $%
\eta _{0}(f)=0$ and $\eta _{0}(f^{2})=1$ we prove that
\begin{equation*}
\sigma _{n}^{2}(F_{n})=\sum_{0\leq p\leq n}~\mathbb{E}\left( \left[
\sum_{0\leq q\leq n}M^{(q-p)_{+}}(f)(X_{q})\right] ^{2}\right)
\end{equation*}%
with the positive part $a_{+}=\max {(a,0)}$ and the convention $M^{0}=Id$,
the identity transition. For $M(x,dy)=\eta _{0}(dy)$, we find that
\begin{equation}
\sigma _{n}^{2}(F_{n})=\sum_{0\leq p\leq n}~\mathbb{E}\left( \left[
\sum_{0\leq q\leq p}f(X_{q})\right] ^{2}\right) ={(n+1)(n+2)}/{2}
\label{additivefunctex}
\end{equation}

We further assume that the Markov transitions $M_n(x_{n-1},dx_n)$ are
absolutely continuous with respect to some measures $\lambda_n(dx_n)$ on $%
E_n $ and we have
\begin{equation*}
(H)\qquad \forall (x_{n-1},x_n)\in \left(E_{n-1}\times E_n\right)\qquad
H_n(x_{n-1},x_n)=\frac{dM_n(x_{n-1},\mbox{\LARGE .})}{d\lambda_n}(x_n)>0
\end{equation*}
In this situation, we have the backward decomposition formula
\begin{equation}  \label{backward}
\mathbb{Q}_n(d(x_0,\ldots,x_n))=\eta_n(dx_n)~ \mathcal{M }%
_{n}(x_{n},d(x_0,\ldots,x_{n-1}))
\end{equation}
with the Markov transitions $\mathcal{M }_n$ defined below
\begin{equation*}
\mathcal{M }_{n}(x_{n},d(x_0,\ldots,x_{n-1})):=\prod_{q=1}^{n}M_{q,%
\eta_{q-1}}(x_q,dx_{q-1})
\end{equation*}
In the above display, $M_{n+1,\eta}$ is the collection of Markov transitions
defined for any $n\geq 0$ and $\eta\in \mathcal{P }(E_n)$ by
\begin{equation}  \label{backwardt}
M_{n+1,\eta}(x,dy)=\frac{1}{\eta\left(G_{n} H_{n+1}(\mbox{\LARGE .},x)\right)%
}~G_{n}(y)~H_{n+1}(y,x)~\eta(dy)
\end{equation}
A detailed proof of this formula and its extended version is provided in
section~\ref{bacw}.

Using the representation in (\ref{backward}), one natural way to approximate
$\mathbb{Q}_n$ is to replace the measures $\eta_n$ with their $N$-particle
approximations $\eta^N_n$. The resulting particle approximation measures, $%
\mathbb{Q}^N_n$, is then
\begin{equation}  \label{backpart}
\mathbb{Q}_n^N(d(x_0,\ldots,x_n)):=\eta^N_n(dx_n)~ \mathcal{M }%
^N_{n}(x_{n},d(x_0,\ldots,x_{n-1}))
\end{equation}
with the random transitions
\begin{equation}  \label{backpartMN}
\mathcal{M }^N_{n}(x_{n},d(x_0,\ldots,x_{n-1})):=\prod_{q=1}^{n}M_{q,%
\eta^N_{q-1}}(x_q,dx_{q-1})
\end{equation}
At this point, it is convenient to recall that for any bounded measurable
function $f_n$ on $E_n$, the measures $\eta_n$ can be written as follows:
\begin{equation}
\eta_n(f_n):=\frac{\gamma_n(f_n)}{\gamma_n(1)}\quad\mbox{\rm with}\quad
\gamma_n(f_n):=\mathbb{E}\left(f_n(X_n)~\prod_{0\leq
p<n}G_p(X_{p})\right)=\eta_n(f_n)~ \prod_{0\leq p<n}\eta_p(G_p)
\label{cvref3}
\end{equation}
The multiplicative formula in the r.h.s. of (\ref{cvref3}) is easily checked
using the fact that $\gamma_{n+1}(1)=\gamma_{n}(G_{n})=\eta_n(G_n)~%
\gamma_n(1)$. Mimicking the above formulae, we set
\begin{equation*}
\Gamma^N_n=\gamma^N_n(1)\times \mathbb{Q}^N_n\quad\mbox{\rm with}\quad
\gamma^N_n(1):=\prod_{0\leq p<n}\eta^N_p(G_p) \quad\mbox{\rm and}\quad
\gamma^N_n(dx)=\gamma^N_n(1)\times \eta^N_n(dx)
\end{equation*}

Notice that the $N$-particle approximation measures $\mathbb{Q}_{n}^{N}$ can
be computed recursively with respect to the time parameter. For instance,
for linear functionals of the form (\ref{additivefunct}), we have
\begin{equation*}
\mathbb{Q}_{n}^{N}(F_{n})=\eta _{n}^{N}(F_{n}^{N})
\end{equation*}%
with a sequence of random functions $F_{n}^{N}$ on $E_{n}$ that can be
computed \textquotedblleft on-the-fly\textquotedblright\ according to the
following recursion
\begin{equation*}
F_{n}^{N}=\sum_{0\leq p\leq n}\left[ M_{n,\eta _{n-1}^{N}}\ldots M_{p+1,\eta
_{p}^{N}}\right] (f_{p})=f_{n}+M_{n,\eta _{n-1}^{N}}(F_{n-1}^{N})
\end{equation*}%
with the initial value $F_{0}^{N}=f_{0}$. In contrast to the genealogical
tree based particle model (\ref{cvarbre}), this new particle algorithm
requires $N^{2}$ computations instead of $N$, in the sense that:
\begin{equation*}
\forall 1\leq j\leq N\qquad F_{n}^{N}(\xi _{n}^{j})=f_{n}(\xi
_{n}^{j})+\sum_{1\leq i\leq N}\frac{G_{n-1}(\xi _{n-1}^{i})H_{n}(\xi
_{n-1}^{i},\xi _{n}^{j})}{\sum_{1\leq i^{\prime }\leq N}G_{n-1}(\xi
_{n-1}^{i^{\prime }})H_{n}(\xi _{n-1}^{i^{\prime }},\xi _{n}^{j})}%
~F_{n-1}^{N}(\xi _{n-1}^{i})
\end{equation*}
This recursion can be straightforwardly extended to the case where we have $%
f_{n}(x_{n-1},x_{n})$ instead of $f_{n}(x_{n})$ in (\ref{additivefunct}) as follows
\begin{equation*}
\forall 1\leq j\leq N\text{ \ \ }F_{n}^{N}(\xi _{n}^{j})=\sum_{1\leq i\leq N}%
\frac{G_{n-1}(\xi _{n-1}^{i})H_{n}(\xi _{n-1}^{i},\xi _{n}^{j})}{\sum_{1\leq
i^{\prime }\leq N}G_{n-1}(\xi _{n-1}^{i^{\prime }})H_{n}(\xi
_{n-1}^{i^{\prime }},\xi _{n}^{j})}~\left( f_{n}(\xi _{n-1}^{i},\xi
_{n}^{j})+F_{n-1}^{N}(\xi _{n-1}^{i})\right)
\end{equation*}%
A very important application of this recursion is to
parameter estimation for non-linear non-Gaussian state-space models.
For instance, it may be used to implement an on-line version of the Expectation-Maximization algorithm as detailed in
\cite[Section 3.2]{kantas2009}.
In a different approach to recursive parameter estimation, an online particle algorithm is presented
in \cite{poyadjis2009} to compute the score for non-linear non-Gaussian
state-space models. In fact, the algorithm of \cite{poyadjis2009} is actually implementing a special
case of the above recursion and may be reinterpreted as an \textquotedblleft
on-the-fly\textquotedblright\ computation of the forward filtering backward
smoothing estimate of an additive functional derived from Fisher's
identity.

The convergence analysis of the $N$-particle measures $\mathbb{Q}_{n}^{N}$
towards their limiting value $\mathbb{Q}_{n}$, as $N\rightarrow \infty $, is
intimately related to the convergence of the flow of particle measures $%
(\eta _{p}^{N})_{0\leq p\leq n}$ towards their limiting measures $(\eta
_{p})_{0\leq p\leq n}$. Several estimates can be easily derived more or less
directly from the convergence analysis of the particle occupation measures $%
\eta _{n}^{N}$ developed in~\cite{fk}, including $\mathbb{L}_{p}$-mean error
bounds and exponential deviation estimates. It is clearly out of the scope
of the present work to review all these consequences. One of the central
objects in this analysis is the local sampling errors $V_{n}^{N}$ induced by
the mean field particle transitions and defined by the following stochastic
perturbation formula
\begin{equation}
\eta _{n}^{N}=\eta _{n-1}^{N}K_{n,\eta _{n-1}^{N}}+\frac{1}{\sqrt{N}}%
~V_{n}^{N}
\end{equation}%
The fluctuation and the deviations of these centered random measures $%
V_{n}^{N}$ can be estimated using non asymptotic Kintchine's type $\mathbb{L}%
_{r}$-inequalities, as well as Hoeffding's or Bernstein's type exponential
deviations~\cite{fk,rio-2009}. We also proved in \cite{dm-2000} that these
random perturbations behave asymptotically as Gaussian random perturbations.
More precisely, for any fixed time horizon $n\geq 0$, the sequence of random
fields $V_{n}^{N}$ converges in law, as the number of particles $N$ tends to
infinity, to a sequence of independent, Gaussian and centered random fields $%
V_{n}$ ; with, for any bounded function $f$ on $E_{n}$, and $n\geq 0$,
\begin{equation}
\mathbb{E}(V_{n}(f)^{2})=\int ~\eta _{n-1}(dx)K_{n,\eta _{n-1}}(x,dy)\left(
f(y)-K_{n,\eta _{n-1}}(f)(x)\right) ^{2}  \label{corr1}
\end{equation}

In section~\ref{secpartic}, we provide some key decompositions expressing
the deviation of the particle measures $(\Gamma _{n}^{N},\mathbb{Q}_{n}^{N})$
around their limiting values $(\Gamma _{n},\mathbb{Q}_{n})$ in terms of
these local random fields models. These decomposition can be used to derive
almost directly some exponential and $\mathbb{L}_{p}$-mean error bounds
using the stochastic analysis developed in~\cite{fk}. We shall use these
functional central limit theorems and some of their variations in various
places in the present article.

\section{Statement of some results}

\label{state} In the present article, we have chosen to concentrate on
functional central limit theorems, as well as on non asymptotic variance
theorems in terms of the time horizon. To describe our results, it is
necessary to introduce the following notation. Let $\beta (M)$ denote the
Dobrushin coefficient of a Markov transition $M$ from a measurable space $E$
into another measurable space $E^{\prime }$ which defined by the following
formula
\begin{equation*}
\beta (M):=\sup {\ \{\mbox{\rm osc}(M(f))\;;\;\;f\in \mbox{\rm Osc}%
_{1}(E^{\prime })\}}
\end{equation*}%
where $\mbox{Osc}_{1}(E^{\prime })$ stands the set of $\mathcal{E^{\prime }}$%
-measurable functions $f$ with oscillation, denoted $\mbox{osc}(f)=\sup {%
\{|f(x)-f(y)|\;;\;x,y\in E}^{\prime }{\}}$, less than or equal to 1. Some
stochastic models discussed in the present article are based on sequences of
random Markov transitions $M^{N}$ that depend on some mean field particle
model with $N$ random particles. In this case, $\beta (M^{N})$ may fail to
be measurable. For this type of models we shall use outer probability
measures to integrate these quantities. For instance, the mean value $%
\mathbb{E}\left( \beta (M^{N})\right) $ is to be understood as the infimum
of the quantities $\mathbb{E}(B^{N})$ where $B^{N}\geq \beta (M^{N})$ are
measurable dominating functions. We also recall that $\gamma _{n}$ satisfy
the linear recursive equation
\begin{equation*}
\gamma _{n}=\gamma _{p}Q_{p,n}\quad \mbox{\rm with}\quad
Q_{p,n}=Q_{p+1}Q_{p+2}\ldots Q_{n}\quad \mbox{\rm and}\quad
Q_{n}(x,dy)=G_{n-1}(x)~M_{n}(x,dy)
\end{equation*}%
for any $0\leq p\leq n$. Using elementary manipulations, we also check that
\begin{equation*}
\Gamma _{n}(F_{n})=\gamma _{p}D_{p,n}(F_{n})
\end{equation*}%
with the bounded integral operators $D_{p,n}$ from $E_{p}$ into $E_{[0,n]}$
defined below
\begin{equation}
D_{p,n}(F_{n})(x_{p}):=\int \mathcal{M}_{p}(x_{p},d(x_{0},\ldots ,x_{p-1}))%
\mathcal{Q}_{p,n}(x_{p},d(x_{p+1},\ldots ,x_{n}))~F_{n}(x_{0},\ldots ,x_{n})
\label{defiDpn}
\end{equation}%
with
\begin{equation*}
\mathcal{Q}_{p,n}(x_{p},d(x_{p+1},\ldots ,x_{n})):=\prod_{p\leq
q<n}Q_{q+1}(x_{q},dx_{q+1})
\end{equation*}

We also let $(G_{p,n},P_{p,n})$ be the pair of potential functions and
Markov transitions defined below
\begin{equation}  \label{DGP}
G_{p,n}=Q_{p,n}(1)/\eta_pQ_{p,n}(1) \quad\mbox{and} \quad P_{p,n}(F_n)={%
D_{p,n}(F_n)}/{D_{p,n}(1)}
\end{equation}
Let the mapping $\Phi_{p,n}: \mathcal{P }(E_p) \rightarrow \mathcal{P }(E_n)
$, $0 \leq p \leq n$, be defined as follows
\begin{equation*}
\Phi_{p,n}(\mu_p)=\frac{\mu_p Q_{p,n}}{\mu_p Q_{p,n}(1)}
\end{equation*}

Our first main result is a functional central limit theorem for the pair of
random fields on $\mathcal{B}(E_{[0,n]})$ defined below
\begin{equation*}
W_{n}^{\Gamma ,N}:=\sqrt{N}\left( \Gamma _{n}^{N}-\Gamma _{n}\right) \quad %
\mbox{\rm and}\quad W_{n}^{\mathbb{Q},N}:=\sqrt{N}~[\mathbb{Q}_{n}^{N}-%
\mathbb{Q}_{n}]
\end{equation*}%
$W_{n}^{\Gamma ,N}$ is centered in the sense that $\mathbb{E}\left(
W_{n}^{\Gamma ,N}(F_{n})\right) =0$ for any $F_{n}\in \mathcal{B}(E_{[0,n]})$%
. The proof of this surprising unbiasedness property can be found in
corollary~\ref{coroo}, in section~\ref{secpartic}.

The first main result of this article is the following multivariate
fluctuation theorem.

\begin{theo}
\label{tclf} We suppose that the following regularity condition is met for
any $n\geq 1$ and for any pair of states $(x,y)\in (E_{n-1},E_{n})$
\begin{equation}
(H^{+})~\quad h_{n}^{-}(y)\leq H_{n}(x,y)\leq h_{n}^{+}(y)\quad \mbox{with}%
\quad ({h_{n}^{+}}/{h_{n}^{-}})\in \mathbb{L}_{4}(\eta _{n})\quad \mbox{and}%
\quad h_{n}^{+}\in \mathbb{L}_{1}(\lambda _{n})  \label{defH+}
\end{equation}%
In this situation, the sequence of random fields $W_{n}^{\Gamma ,N}$, resp. $%
W_{n}^{\mathbb{Q},N}$, converge in law, as $N\rightarrow \infty $, to the
centered Gaussian fields $W_{n}^{\Gamma }$, resp. $W_{n}^{\mathbb{Q}}$,
defined for any $F_{n}\in \mathcal{B}(E_{[0,n]})$ by
\begin{eqnarray*}
W_{n}^{\Gamma }(F_{n}) &=&\sum_{p=0}^{n}\gamma _{p}(1)~V_{p}\left(
D_{p,n}(F_{n})\right) \\
W_{n}^{\mathbb{Q}}(F_{n}) &=&\sum_{p=0}^{n}V_{p}\left( G_{p,n}~P_{p,n}(F_{n}-%
\mathbb{Q}_{n}(F_{n}))\right)
\end{eqnarray*}
\end{theo}

The second main result of the article is the following non asymptotic
theorem.

\begin{theo}
\label{nonasymptheo} For any $r\geq 1$, $n\geq 0$, $F_n\in \mathcal{B }%
(E_{[0,n]})$ s.t. $\|F_n\|\leq 1$
\begin{equation}  \label{LLr}
\sqrt{N}~\mathbb{E}\left(\left| [\mathbb{Q}^N_n-\mathbb{Q}_n](F_n)
\right|^r\right)^{\frac{1}{r}}\leq a_r~ \sum_{0\leq p\leq
n}~b_{p,n}^2~c^N_{p,n}
\end{equation}
for some finite constants $a_r<\infty$ whose values only depend on the
parameter $r$, and a pair of constants $(b_{p,n},c^N_{p,n})$ such that
\begin{equation*}
b_{p,n}\leq \sup_{x,y}{({Q_{p,n}(1)(x)}/{Q_{p,n}(1)(y)})}\quad\mbox{and}%
\quad c^N_{p,n}\leq \mathbb{E}\left(\beta(P^N_{p,n})\right)
\end{equation*}
In the above display, $P^N_{p,n}$ stands for the random Markov transitions
defined as $P_{p,n}$ by replacing in (\ref{defiDpn}) and (\ref{DGP}) the
transitions $\mathcal{M }_p$ by $\mathcal{M }^N_p$. For linear functionals
of the form (\ref{additivefunct}), with $f_n\in\mbox{\rm Osc}_{1}(E_n)$, the
constant $c^N_{p,n}$ in (\ref{LLr}) can be chosen so that
\begin{equation}  \label{estitheo}
c^N_{p,n}\leq \sum_{0\leq q<p}\beta\left(M_{p,\eta^N_{p-1}}\ldots
M_{q+1,\eta^N_q}\right) +\sum_{p\leq q\leq n}~b_{q,n}^2~\beta(S_{p,q})
\end{equation}
with the Markov transitions $S_{p,q}$ from $E_p$ into $E_q$ defined for any
function $f\in \mathcal{B }(E_q)$ by the following formula $%
S_{p,q}(f)=Q_{p,q}(f)/Q_{p,q}(1)$.
\end{theo}

We emphasize that the $\mathbb{L}_{r}$-mean error bounds described in the
above theorem enter the stability properties of the semigroups $S_{p,q}$ and
the one associated with the backward Markov transitions $M_{n+1,\eta
_{n}^{N}}$. In several instances, the term in the r.h.s. of (\ref{estitheo})
can be uniformly bounded with respect to the time horizon. For instance, in
the toy example we discussed in (\ref{additivefunctex}), we have the
estimates
\begin{equation*}
b_{p,n}=1\quad \mbox{\rm and}\quad c_{p,n}^{N}\leq 1~\Longrightarrow ~\sqrt{N%
}~\mathbb{E}\left( \left\vert [\mathbb{Q}_{n}^{N}-\mathbb{Q}%
_{n}](F_{n})\right\vert ^{r}\right) ^{\frac{1}{r}}\leq a_{r}~(n+1)
\end{equation*}%
In more general situations, these estimates are related to the stability
properties of the Feynman-Kac semigroup. To simplify the presentation, let
us suppose that the pair of potential-transitions $(G_{n},M_{n})$ are time
homogeneous $(G_{n},H_{n},M_{n})=(G,H,M)$ and chosen so that the following
regularity condition is satisfied
\begin{equation}
(M)_{m}\qquad \forall (x,x^{\prime 2}\quad G(x)\leq \delta ~G(x^{\prime
})\quad \mbox{\rm and}\quad M^{m}(x,dy)\leq \rho ~M^{m}(x^{\prime },dy)
\label{condHh}
\end{equation}%
for some $m\geq 1$ and some parameters $(\delta ,\rho )\in \lbrack 1,\infty
)^{2}$. Under this rather strong condition, we have
\begin{equation*}
b_{p,n}\leq \rho \delta ^{m}\quad \mbox{\rm and}\quad \beta (S_{p,q})\leq
\left( 1-\rho ^{-2}\delta ^{-m}\right) ^{\lfloor (q-p)/m\rfloor }
\end{equation*}%
See for instance corollary 4.3.3. in~\cite{fk} and the more recent article~%
\cite{cdg}. On the other hand, let us suppose that
\begin{equation*}
\inf_{x,y,y^{\prime }}{(H(x,y)/H(x,y^{\prime }))}=\alpha (h)>0
\end{equation*}%
In this case, we have
\begin{equation*}
M_{n,\eta }(x,dy)\leq \alpha (h)^{-2}~M_{n,\eta }(x^{\prime
},dy)\Longrightarrow \beta \left( M_{p,\eta _{p-1}^{N}}\ldots M_{q+1,\eta
_{q}^{N}}\right) \leq \left( 1-\alpha (h)^{2}\right) ^{p-q}
\end{equation*}%
For linear functional models of the form (\ref{additivefunct}) associated
with functions $f_{n}\in \mbox{\rm Osc}_{1}(E_{n})$, it is now readily
checked that
\begin{equation}
\sqrt{N}~\mathbb{E}\left( \left\vert [\mathbb{Q}_{n}^{N}-\mathbb{Q}%
_{n}](F_{n})\right\vert ^{r}\right) ^{\frac{1}{r}}\leq a_{r}~b~(n+1)
\label{fineintro}
\end{equation}%
for some finite constant $b<\infty $ whose values do not depend on the time
parameter $n$. With some information on the constants $a_{r}$, these $%
\mathbb{L}_{r}$-mean error bounds can turned to uniform exponential
estimates w.r.t. the time parameter for normalized additive functionals of
the following form
\begin{equation*}
\overline{F}_{n}(x_{0},\ldots ,x_{n}):=\frac{1}{n+1}\sum_{0\leq p\leq
n}f_{p}(x_{p})
\end{equation*}%
To be more precise, by lemma 7.3.3 in~\cite{fk}, the collection of constants
$a_{r}$ in (\ref{fineintro}) can be chosen so that
\begin{equation}
a_{2r}^{2r}\leq (2r)!~2^{-r}/r!\quad \mbox{and}\quad a_{2r+1}^{2r+1}\leq
(2r+1)!~2^{-r}/r!  \label{amdef}
\end{equation}%
In this situation, it is easily checked that for any $\epsilon >0$, and $%
N\geq 1$, we have the following uniform Gaussian concentration estimates:
\begin{equation*}
{\frac{1}{N}\log {\sup_{n\geq 0}{\mathbb{P}\left( \left\vert [\mathbb{Q}%
_{n}^{N}-\mathbb{Q}_{n}](\overline{F}_{n})\right\vert \geq \frac{b}{\sqrt{N}}%
+\epsilon \right) }}}\leq -{\epsilon ^{2}}/{(2b^{2})}
\end{equation*}%
This result is a direct consequence of the fact that for any non negative
random variable $U$
\begin{equation*}
\left( \forall r\geq 1\quad \mathbb{E}\left( U^{r}\right) ^{\frac{1}{r}}\leq
a_{r}~b\right) \Rightarrow \log {\mathbb{P}\left( U\geq b+\epsilon \right) }%
\leq -{\epsilon ^{2}}/{(2b^{2})}
\end{equation*}%
To check this claim, we develop the exponential to prove that
\begin{equation*}
\log {\mathbb{E}\left( e^{tU}\right) }\overset{\forall t\geq 0}{\leq }bt+%
\frac{(bt)^{2}}{2}\Rightarrow \log {\mathbb{P}\left( U\geq b+\epsilon
\right) }\leq -\sup_{t\geq 0}{\left( \epsilon t-\frac{(bt)^{2}}{2}\right) }
\end{equation*}

\section{A backward Markov chain formulation}

\label{bacw}

This section is mainly concerned with the proof of the backward
decomposition formula (\ref{backward}). Before proceeding, we recall that
the measures $(\gamma _{n},\eta _{n})$ satisfy the non linear equations
\begin{equation*}
\gamma _{n}=\gamma _{n-1}Q_{n}\quad \mbox{\rm and}\quad \eta _{n+1}:=\Phi
_{n+1}(\eta _{n}):=\Psi _{G_{n}}(\eta _{n})M_{n+1}
\end{equation*}%
and their semigroups are given by
\begin{equation*}
\gamma _{n}=\gamma _{p}Q_{p,n}\quad \mbox{\rm and}\quad \eta _{n}(f_{n}):={%
\eta _{p}Q_{p,n}(f_{n})}/{\eta _{p}Q_{p,n}(1)}
\end{equation*}%
for any function $f_{n}\in \mathcal{B}(E_{n})$. In this connection, we also
mention that the semigroup of the pair of measures $(\Gamma _{n},\mathbb{Q}%
_{n})$ defined in (\ref{defQn}) for any $0\leq p\leq n$ and any $F_{n}\in
\mathcal{B}(E_{[0,n]})$, we have
\begin{equation}
\Gamma _{n}(F_{n})=\gamma _{p}D_{p,n}(F_{n})\quad \mbox{\rm and}\quad
\mathbb{Q}_{n}(F_{n})={\eta _{p}D_{p,n}(F_{n})}/{\eta _{p}D_{p,n}(1)}
\label{sgQ}
\end{equation}%
These formulae are a direct consequence of the following observation
\begin{equation*}
\eta _{p}D_{p,n}(F_{n})=\int \mathbb{Q}_{p}(d(x_{0},\ldots ,x_{p}))~\mathcal{%
Q}_{p,n}(x_{p},d(x_{p+1},\ldots ,x_{n}))F_{n}(x_{0},\ldots ,x_{n})
\end{equation*}

\begin{lem}
\label{lemback}For any $0\leq p<n$, we have%
\begin{equation}
\gamma _{p}(dx_{p})~\mathcal{Q}_{p,n}(x_{p},d(x_{p+1},\ldots ,x_{n}))=\gamma
_{n}(dx_{n})~\mathcal{M}_{n,p}(x_{n},d(x_{p},\ldots ,x_{n-1}))
\label{backg2}
\end{equation}%
with
\begin{equation*}
\mathcal{M}_{n,p}(x_{n},d(x_{p},\ldots ,x_{n-1})):=\prod_{p\leq
q<n}M_{q+1,\eta _{q}}(x_{q+1},dx_{q})
\end{equation*}%
In particular, for any time $n\geq 0$, the Feynman-Kac path measures $%
\mathbb{Q}_{n}$ defined in (\ref{defQn}) can be expressed in terms of the
sequence of marginal measures $(\eta _{p})_{0\leq p\leq n}$, with the
following backward Markov chain formulation
\begin{equation}
\mathbb{Q}_{n}(d(x_{0},\ldots ,x_{n}))=\eta _{n}(dx_{n})~\mathcal{M}%
_{n,0}(x_{n},d(x_{0},\ldots ,x_{n-1})):  \label{back2}
\end{equation}
\end{lem}

Before entering into the details of the proof of this lemma, we mention that
(\ref{back2}) holds true for any well defined Markov transition $M_{n+1,\eta
_{n}}(y,dx)$ from $E_{n}$ into $E_{n+1}$ satisfying the local backward
equation
\begin{equation*}
\Psi _{G_{n}}(\eta _{n})(dx)~M_{n+1}(x,dy)=\Phi _{n+1}(\eta
_{n})(dy)~M_{n+1,\eta _{n}}(y,dx)
\end{equation*}%
or equivalently
\begin{equation}
\eta _{n}(dx)~Q_{n+1}(x,dy)=(\eta _{n}Q_{n+1})(dy)~M_{n+1,\eta _{n}}(y,dx)
\label{dualback}
\end{equation}%
In other words, we have the duality formula
\begin{equation}
\Psi _{G_{n}}(\eta _{n})\left( f~M_{n+1}(g)\right) =\Phi _{n+1}(\eta
_{n})\left( g~M_{n+1,\eta _{n}}(f)\right)  \label{dual}
\end{equation}

Also notice that for any pair of measures $\mu ,\nu $ on $E_{n}$ s.t. $\mu
\ll \nu $, we have $\mu M_{n+1}\ll \nu M_{n+1}$. Indeed, if we have $\nu
M_{n+1}(A)=0$, the function $M_{n+1}(1_{A})$ is null $\nu $-almost
everywhere, and therefore $\mu $-almost everywhere from which we conclude
that $\mu M_{n+1}(A)=0$. For any bounded measurable function $g$ on $E_{n}$
we set
\begin{equation*}
\Psi _{G_{n}}^{g}(\eta _{n})(dx)=\Psi _{G_{n}}(\eta _{n})(dx)~g(x)\ll \Psi
_{G_{n}}(\eta _{n})(dx)
\end{equation*}%
>From the previous discussion, we have $\Psi _{G_{n}}^{g}(\eta
_{n})M_{n+1}\ll \Psi _{G_{n}}(\eta _{n})M_{n+1}$ and it is easily checked
that
\begin{equation*}
M_{n+1,\eta _{n}}(g)(y)=\frac{d\Psi _{G_{n}}^{g}(\eta _{n})M_{n+1}}{d\Psi
_{G_{n}}(\eta _{n})M_{n+1}}(y)
\end{equation*}%
is a well defined Markov transition from $E_{n+1}$ into $E_{n}$ satisfying
the desired backward equation. These manipulations are rather classical in
the literature on Markov chains (see for instance~\cite{revuz}, and
references therein). Under the regularity condition $(H)$ the above
transition is explicitly given by the formula (\ref{backwardt}).

Now, we come to the proof of lemma~\ref{lemback}.

\textbf{Proof of lemma~\ref{lemback}:}

We prove (\ref{backg2}) using a backward induction on the parameter $p$. By (%
\ref{dualback}), the formula is clearly true for $p=(n-1)$. Suppose the
result has been proved at rank $p$. Since we have
\begin{equation*}
\begin{array}{l}
\gamma _{p-1}(dx_{p-1})~\mathcal{Q}_{p-1,n}(x_{p-1},d(x_{p},\ldots ,x_{n}))
\\
\\
=\gamma _{p-1}(dx_{p-1})~Q_{p}(x_{p-1},dx_{p})~\mathcal{Q}%
_{p,n}(x_{p},d(x_{p+1},\ldots ,x_{n}))%
\end{array}%
\end{equation*}%
and
\begin{equation*}
\gamma _{p-1}(dx_{p-1})~Q_{p}(x_{p-1},dx_{p})=\gamma _{p}(dx_{p})~M_{p,\eta
_{p-1}}(x_{p},dx_{p-1})
\end{equation*}%
Using the backward induction we conclude that the desired formula is also
met at rank $(p-1)$. The second assertion is a direct consequence of (\ref%
{backg2}). The end of the proof of the lemma is now completed. \hfill %
\hbox{\vrule height 5pt width 5pt depth 0pt}\medskip \newline

We end this section with some properties of backward Markov transitions
associated with a given initial probability measure that may differ from the
one associated with the Feynman-Kac measures. These mathematical objects
appear in a natural way in the analysis of the $N$-particle approximation
transitions $\mathcal{M}_{n}^{N}$ introduced in (\ref{backpartMN}).

\begin{defi}
\label{defiMQD} For any $0\leq p\leq n$ and any probability measure $\eta\in
\mathcal{P }(E_{p})$, we denote by $\mathcal{M }_{n+1,p,\eta}$ the Markov
transition from $E_{n+1}$ into $E_{[p,n]}=(E_p\times\ldots\times E_n)$
defined by
\begin{equation*}
\mathcal{M }_{n+1,p,\eta}\left(x_{n+1},d(x_p,\ldots,x_n)\right)
=\prod_{p\leq q\leq n} M_{q+1,\Phi_{p,q}(\eta)}(x_{q+1},dx_q)
\end{equation*}
\end{defi}

Notice that this definition is consistent with the definition of the Markov
transitions $\mathcal{M }_{p,n}$ introduced in lemma~\ref{lemback}:
\begin{equation*}
\mathcal{M }_{n+1,p,\eta_p}\left(x_{n+1},d(x_p,\ldots,x_n)\right)= \mathcal{%
M }_{n+1,p}\left(x_{n+1},d(x_p,\ldots,x_n)\right)
\end{equation*}
Also observe that $\mathcal{M }_{n+1,p,\eta}$ can alternatively be defined
by the pair of recursions
\begin{equation}  \label{recursions}
\begin{array}{l}
\mathcal{M }_{n+1,p,\eta}\left(x_{n+1},d(x_p,\ldots,x_n)\right) \\
\\
= \mathcal{M }_{n+1,p+1,\Phi_{p+1}(\eta)}\left(x_{n+1},d(x_{p+1},\ldots,x_n)%
\right)\times M_{p+1,\eta}(x_{p+1},dx_p) \\
\\
=M_{n+1,\Phi_{p,n}(\eta)}(x_{n+1},dx_n)~ \mathcal{M }_{n,p,\eta}%
\left(x_{n},d(x_p,\ldots,x_{n-1})\right)%
\end{array}%
\end{equation}
The proof of the following lemma follows the same lines of arguments as the
ones used in the proof of lemma~\ref{lemback}. For the convenience of the
reader, the details of this proof are postponed to the appendix.

\begin{lem}
\label{lemMeta} For any $0\leq p<n$ and any probability measure $\eta \in
\mathcal{P}(E_{p})$, we have
\begin{equation*}
\eta Q_{p,n}(dx_{n})~\mathcal{M}_{n,p,\eta }(x_{n},d(x_{p},\ldots
,x_{n-1}))=\eta (dx_{p})~\mathcal{Q}_{p,n}(x_{p},d(x_{p+1},\ldots ,x_{n}))
\end{equation*}%
In other words, we have
\begin{equation}
\begin{array}{l}
\mathcal{M}_{n,p,\eta }(x_{n},d(x_{p},\ldots ,x_{n-1})) \\
\\
=\displaystyle\frac{(\eta \times \mathcal{Q}_{p,n-1})(d(x_{p},\ldots
,x_{n-1}))G_{n-1}(x_{n-1})~H_{n}(x_{n-1},x_{n})}{(\eta Q_{p,n-1})\left(
G_{n-1}~H_{n}(\mbox{\LARGE .},x_{n})\right) }%
\end{array}
\label{lemformu}
\end{equation}%
with the measure $(\eta \times \mathcal{Q}_{p,n-1})$ defined below
\begin{equation*}
(\eta \times \mathcal{Q}_{p,n-1})(d(x_{p},\ldots ,x_{n-1})):=\eta (dx_{p})~%
\mathcal{Q}_{p,n-1}(x_{p},d(x_{p+1},\ldots ,x_{n-1}))
\end{equation*}
\end{lem}

\section{Particle approximation models}

\label{secpartic}

We provide in this section some preliminary results on the convergence of
the $N$-particle measures $(\Gamma^N_n,\mathbb{Q}^N_n)$ to their limiting
values $(\Gamma_n,\mathbb{Q}_n)$, as $N\rightarrow\infty$. Most of the
forthcoming analysis is developed in terms of the following integral
operators.

\begin{defi}
For any $0\leq p\leq n$, we let $D^N_{p,n}$ be the $\mathcal{F }^N_{p-1}$%
-measurable integral operators from $\mathcal{B }(E_{[0,n]})$ into $\mathcal{%
B }(E_p)$ defined below
\begin{equation*}
D^N_{p,n}(F_n)(x_p):=\int \mathcal{M }_p^N(x_{p},d(x_0,\ldots,x_{p-1}))
\mathcal{Q }_{p,n}(x_p,d(x_{p+1},\ldots,x_{n})) F_n(x_0,\ldots,x_n)
\end{equation*}
with the conventions $D^N_{0,n}= \mathcal{Q }_{0,n}$, and resp. $D^N_{n,n}=
\mathcal{M }_n^N$, for $p=0$, and resp. $p=n$
\end{defi}

The main result of this section is the following theorem.

\begin{theo}
\label{theo2} For any $0\leq p\leq n$, and any function $F_{n}$ on the path
space $E_{[0,n]}$, we have
\begin{equation*}
\mathbb{E}\left( \Gamma _{n}^{N}(F_{n})\left\vert ~\mathcal{F}%
_{p}^{N}\right. \right) =\gamma _{p}^{N}\left( D_{p,n}^{N}(F_{n})\right)
\quad \mbox{and}\quad W_{n}^{\Gamma ,N}(F_{n})=\sum_{p=0}^{n}\gamma
_{p}^{N}(1)~V_{p}^{N}\left( D_{p,n}^{N}(F_{n})\right)
\end{equation*}
\end{theo}

\textbf{Proof of theorem~\ref{theo2}:}

To prove the first assertion, we use a backward induction on the parameter $%
p $. For $p=n$, the result is immediate since we have
\begin{equation*}
\Gamma _{n}^{N}(F_{n})=\gamma _{n}^{N}(1)~\eta _{n}^{N}\left(
D_{n,n}^{N}(F_{n})\right)
\end{equation*}%
We suppose that the formula is valid at a given rank $p\leq n$. In this
situation, we have
\begin{eqnarray}
\mathbb{E}\left( \Gamma _{n}^{N}(F_{n})\left\vert ~\mathcal{F}%
_{p-1}^{N}\right. \right) &=&\gamma _{p}^{N}(1)~\mathbb{E}\left( \eta
_{p}^{N}\left( D_{p,n}^{N}(F_{n})\right) \left\vert ~\mathcal{F}%
_{p-1}^{N}\right. \right)  \notag \\
\ &=&\gamma _{p-1}^{N}(1)~\int \eta _{p-1}^{N}(G_{p-1}H_{p}(\mbox{\LARGE .}%
,x_{p}))~\lambda _{p}(dx_{p})~D_{p,n}^{N}(F_{n})(x_{p})  \label{refrhs}
\end{eqnarray}%
Using the fact that
\begin{equation*}
\gamma _{p-1}^{N}(1)~\eta _{p-1}^{N}(G_{p-1}H_{p}(\mbox{\LARGE .}%
,x_{p}))~\lambda _{p}(dx_{p})~M_{p,\eta _{p-1}^{N}}(x_{p},dx_{p-1})=\gamma
_{p-1}^{N}(dx_{p-1})Q_{p}(x_{p-1},dx_{p})
\end{equation*}%
we conclude that the r.h.s. term in (\ref{refrhs}) takes the form
\begin{equation*}
\begin{array}{l}
\int \gamma _{p-1}^{N}(dx_{p-1})\mathcal{M}_{p-1}^{N}(x_{p-1},d(x_{0},\ldots
,x_{p-2}))\mathcal{Q}_{p-1,n}(x_{p-1},d(x_{p},\ldots
,x_{n}))~F_{n}(x_{0},\ldots ,x_{n}) \\
\\
=\gamma _{p-1}^{N}\left( D_{p-1,n}^{N}(F_{n})\right)%
\end{array}%
\end{equation*}%
This ends the proof of the first assertion. The proof of the second
assertion is based on the following decomposition
\begin{eqnarray*}
\left( \Gamma _{n}^{N}-\Gamma _{n}\right) (F_{n}) &=&\sum_{p=0}^{n}\left[
\mathbb{E}\left( \Gamma _{n}^{N}(F_{n})\left\vert ~\mathcal{F}%
_{p}^{N}\right. \right) -\mathbb{E}\left( \Gamma _{n}^{N}(F_{n})\left\vert ~%
\mathcal{F}_{p-1}^{N}\right. \right) \right] \\
&=&\sum_{p=0}^{n}\gamma _{p}^{N}(1)~\left( \eta _{p}^{N}\left(
D_{p,n}^{N}(F_{n})\right) -\frac{1}{\eta _{p-1}^{N}(G_{p-1})}~\eta
_{p-1}^{N}\left( D_{p-1,n}^{N}(F_{n})\right) \right)
\end{eqnarray*}%
where $\mathcal{F}_{-1}^{N}$ is the trivial sigma field. By definition of
the random fields $V_{p}^{N}$, it remains to prove that
\begin{equation*}
\eta _{p-1}^{N}\left( D_{p-1,n}^{N}(F_{n})\right) =(\eta
_{p-1}^{N}Q_{p})\left( D_{p,n}^{N}(F_{n})\right)
\end{equation*}%
To check this formula, we use the decomposition
\begin{equation}
\begin{array}{l}
\eta _{p-1}^{N}(dx_{p-1})~\mathcal{M}_{p-1}^{N}(x_{p-1},d(x_{0},\ldots
,x_{p-2}))~\mathcal{Q}_{p-1,n}(x_{p-1},d(x_{p},\ldots ,x_{n})) \\
\\
=\eta _{p-1}^{N}(dx_{p-1})Q_{p}(x_{p-1},dx_{p})\mathcal{M}%
_{p-1}^{N}(x_{p-1},d(x_{0},\ldots ,x_{p-2}))~\mathcal{Q}%
_{p,n}(x_{p},d(x_{p+1},\ldots ,x_{n}))%
\end{array}
\label{refeqq}
\end{equation}

Using the fact that
\begin{eqnarray*}
\eta_{p-1}^N(dx_{p-1}) Q_p(x_{p-1},dx_{p})&=&(\eta^N_{p-1}Q_p)(dx_p)~
M_{p,\eta^N_{p-1}}(x_p,dx_{p-1})~
\end{eqnarray*}
we conclude that the term in the r.h.s. of (\ref{refeqq}) is equal to
\begin{equation*}
\begin{array}{l}
(\eta^N_{p-1}Q_p)(dx_p)~ \mathcal{M }^N_{p}(x_{p},d(x_0,\ldots,x_{p-1}))~
\mathcal{Q }_{p,n}(x_{p},d(x_{p+1},\ldots,x_n))%
\end{array}%
\end{equation*}
This ends the proof of the theorem. \hfill%
\hbox{\vrule height 5pt width 5pt
depth 0pt}\medskip \newline

Several consequences of theorem~\ref{theo2} are now emphasized. On the one
hand, using the fact that the random fields $V_{n}^{N}$ are centered given $%
\mathcal{F}_{n-1}^{N}$, we find that
\begin{equation*}
\mathbb{E}\left( \Gamma _{n}^{N}(F_{n})\right) =\Gamma _{n}(F_{n})
\end{equation*}%
On the other hand, using the fact that
\begin{equation*}
\frac{\gamma _{p}(1)}{\gamma _{n}(1)}=\frac{\gamma _{p}(1)}{\gamma
_{p}Q_{p,n}(1)}=\frac{1}{\eta _{p}Q_{p,n}(1)}
\end{equation*}%
we prove the following decomposition
\begin{equation}
\overline{W}_{n}^{\Gamma ,N}(F_{n})=\sqrt{N}\left( \overline{\gamma }%
_{n}^{N}(1)~\mathbb{Q}_{n}^{N}-\mathbb{Q}_{n}\right) (F_{n})=\sum_{p=0}^{n}%
\overline{\gamma }_{p}^{N}(1)~V_{p}^{N}\left( \overline{D}%
_{p,n}^{N}(F_{n})\right)  \label{overlW}
\end{equation}%
with the pair of parameters $\left( \overline{\gamma }_{n}^{N}(1),\overline{D%
}_{p,n}^{N}\right) $ defined below
\begin{equation}
\overline{\gamma }_{n}^{N}(1):=\frac{\gamma _{n}^{N}(1)}{\gamma _{n}(1)}%
\quad \mbox{\rm
and}\quad \overline{D}_{p,n}^{N}(F_{n})=\frac{D_{p,n}^{N}(F_{n})}{\eta
_{p}Q_{p,n}(1)}  \label{overlDN}
\end{equation}%
Using again the fact that the random fields $V_{n}^{N}$ are centered given $%
\mathcal{F}_{n-1}^{N}$, we have
\begin{equation*}
\mathbb{E}\left( \overline{W}_{n}^{\Gamma ,N}(F_{n})^{2}\right)
=\sum_{p=0}^{n}\mathbb{E}\left( \overline{\gamma }_{p}^{N}(1)^{2}~\mathbb{E}%
\left[ V_{p}^{N}\left( \overline{D}_{p,n}^{N}(F_{n})\right) ^{2}~\left\vert ~%
\mathcal{F}_{p-1}^{N}\right. \right] \right)
\end{equation*}%
Using the estimates
\begin{eqnarray}
\Vert D_{p,n}^{N}(F_{n})\Vert &\leq &\Vert Q_{p,n}(1)\Vert ~\Vert F_{n}\Vert
\notag \\
\Vert \overline{D}_{p,n}^{N}(F_{n})\Vert &\leq &\Vert \overline{Q}%
_{p,n}(1)\Vert ~\Vert F_{n}\Vert \quad \mbox{\rm with}\quad \overline{Q}%
_{p,n}(1)=\frac{Q_{p,n}(1)}{\eta _{p}Q_{p,n}(1)}  \label{defqover}
\end{eqnarray}%
we prove the non asymptotic variance estimate
\begin{equation*}
\mathbb{E}\left( \overline{W}_{n}^{\Gamma ,N}(F_{n})^{2}\right) \leq
\sum_{p=0}^{n}\mathbb{E}\left( \overline{\gamma }_{p}^{N}(1)^{2}\right)
~\Vert \overline{Q}_{p,n}(1)\Vert ^{2}=\sum_{p=0}^{n}\left[ 1+\mathbb{E}%
\left( [\overline{\gamma }_{p}^{N}(1)-1]^{2}\right) \right] ~\Vert \overline{%
Q}_{p,n}(1)\Vert ^{2}
\end{equation*}%
for any function $F_{n}$ such that $\Vert F_{n}\Vert \leq 1$. On the other
hand, using the decomposition
\begin{equation*}
\left( \overline{\gamma }_{n}^{N}(1)~\mathbb{Q}_{n}^{N}-\mathbb{Q}%
_{n}\right) =\left[ \overline{\gamma }_{n}^{N}(1)-1\right] ~\mathbb{Q}%
_{n}^{N}+\left( \mathbb{Q}_{n}^{N}-\mathbb{Q}_{n}\right)
\end{equation*}%
we prove that
\begin{equation*}
\mathbb{E}\left( \left[ \mathbb{Q}_{n}^{N}(F_{n})-\mathbb{Q}_{n}(F_{n})%
\right] ^{2}\right) ^{1/2}\leq \frac{1}{\sqrt{N}}~\mathbb{E}\left(
W_{n}^{\Gamma }(F_{n})^{2}\right) ^{1/2}+\mathbb{E}\left( \left[ \overline{%
\gamma }_{n}^{N}(1)-1\right] ^{2}\right) ^{1/2}
\end{equation*}%
Some interesting bias estimates can also be obtained using the fact that
\begin{equation*}
\mathbb{E}\left( \mathbb{Q}_{n}^{N}(F_{n})\right) -\mathbb{Q}_{n}(F_{n})=%
\mathbb{E}\left( \left[ 1-\overline{\gamma }_{n}^{N}(1)\right] ~\left[
\mathbb{Q}_{n}^{N}(F_{n})-\mathbb{Q}_{n}(F_{n})\right] \right)
\end{equation*}%
and the following easily proved upper bound
\begin{equation*}
\left\vert \mathbb{E}\left( \mathbb{Q}_{n}^{N}(F_{n})\right) -\mathbb{Q}%
_{n}(F_{n})\right\vert \leq \mathbb{E}\left( \left[ 1-\overline{\gamma }%
_{n}^{N}(1)\right] ^{2}\right) ^{1/2}\mathbb{E}\left( \left[ \mathbb{Q}%
_{n}^{N}(F_{n})-\mathbb{Q}_{n}(F_{n})\right] ^{2}\right) ^{1/2}
\end{equation*}

Under the regularity condition $(M)_{m}$ stated in (\ref{condHh}), we proved
in a recent article~\cite{cdg}, that for any $n\geq p\geq 0$, and any $%
N>(n+1)\rho \delta ^{m}$ we have
\begin{equation*}
\Vert \overline{Q}_{p,n}(1)\Vert \leq \delta ^{m}\rho \quad \mbox{\rm and}%
\quad N~\mathbb{E}\left[ \left( \overline{\gamma }_{n}^{N}(1)-1\right) ^{2}%
\right] \leq 4~(n+1)~\rho ~\delta ^{m}
\end{equation*}%
>From these estimates, we readily prove the following corollary.

\begin{cor}
\label{coroo} Assume that condition $(M)_m$ is satisfied for some parameters
$(m,\delta,\rho)$. In this situation, for any $n\geq p\geq 0$, any $F_n$
such that $\|F_n\|\leq 1$, and any $N> (n+1)\rho\delta^{m}$ we have
\begin{equation*}
\mathbb{E}\left(\overline{W}^{\Gamma,N}_n(F_n)\right)=0\quad\mbox{ and}\quad
\mathbb{E}\left(\overline{W}^{\Gamma,N}_n(F_n)^2\right)\leq (
\delta^m\rho)^2 (n+1) \left( 1+\frac{2}{N}~\rho \delta^{m} (n+2) \right)
\end{equation*}
In addition, we have
\begin{eqnarray*}
N~\mathbb{E}\left(\left[\mathbb{Q}^N_n(F_n)-\mathbb{Q}_n(F_n)\right]%
^2\right) &\leq & 2 (n+1) \rho\delta^m\left( 4+\rho \delta^{m}\left[1 +\frac{%
2}{N} (n+2)\right]\right)
\end{eqnarray*}
and the bias estimate
\begin{equation*}
N\left|\mathbb{E}\left(\mathbb{Q}^N_n(F_n)\right)-\mathbb{Q}%
_n(F_n)\right|\leq 2\sqrt{2}~(n+1)\rho\delta^{m} \left( 4+\rho\delta^{m}%
\left[1 +\frac{2}{N} (n+2)\right]\right)^{1/2}
\end{equation*}
\end{cor}

\section{Fluctuation properties}

\label{fluctsec} This section is mainly concerned with the proof of theorem~%
\ref{tclf}. Unless otherwise is stated, in the further developments of this
section, we assume that the regularity condition $(H^{+})$ presented in (\ref%
{defH+}) is satisfied for some collection of functions $%
(h_{n}^{-},h_{n}^{+}) $. Our first step to establish theorem~\ref{tclf} is
the fluctuation analysis of the $N$-particle measures $(\Gamma _{n}^{N},%
\mathbb{Q}_{n}^{N})$ given in proposition \ref{propDN} whose proof relies on
the following technical lemma.

\begin{lem}
\label{lemLN}
\begin{equation*}
\begin{array}{l}
\mathcal{M}_{n}^{N}(x_{n},d(x_{0},\ldots ,x_{n-1}))-\mathcal{M}%
_{n}(x_{n},d(x_{0},\ldots ,x_{n-1})) \\
\\
=\displaystyle\sum_{0\leq p\leq n}\left[ \mathcal{M}_{n,p,\eta _{p}^{N}}-%
\mathcal{M}_{n,p,\Phi _{p}\left( \eta _{p-1}^{N}\right) }\right] \left(
x_{n},d(x_{p},\ldots ,x_{n-1})\right) ~\mathcal{M}_{p}^{N}(x_{p},d(x_{0},%
\ldots ,x_{p-1}))%
\end{array}%
\end{equation*}
\end{lem}

The proof of this lemma follows elementary but rather tedious calculations;
thus it is postponed to the appendix. We now state proposition \ref{propDN}.

\begin{prop}
\label{propDN} For any $N\geq 1$, $0\leq p\leq n$, $x_p\in E_p$, $m\geq 1$,
and $F_n\in \mathcal{B }(E_{[0,n]})$ such that $\|F_n\|\leq 1$, we have
\begin{equation}  \label{keyest}
\sqrt{N}~\mathbb{E}\left(\left|D_{p,n}^N(F_n)-D_{p,n}(F_n)(x_p)\right|^{m}
\right)^{\frac{1}{m}} \leq a(m)~b(n)~\left(\frac{h^+_p}{h^-_p}(x_p)\right)^2
\end{equation}
for some finite constants $a(m)<\infty$, resp. $b(n)<\infty$, whose values
only depend on the parameters $m$, resp. on the time horizon $n$.
\end{prop}

\noindent\mbox{\bf Proof:}\newline

Using lemma~\ref{lemLN}, we find that
\begin{equation*}
D_{p,n}^{N}(F_{n})-D_{p,n}(F_{n})=\displaystyle\sum_{0\leq q\leq p}\left[
\mathcal{M}_{p,q,\eta _{q}^{N}}-\mathcal{M}_{p,q,\Phi _{q}\left( \eta
_{q-1}^{N}\right) }\right] \left( T_{p,q,n}^{N}(F_{n})\right)
\end{equation*}%
with the random function $T_{p,q,n}^{N}(F_{n})$ defined below
\begin{equation*}
\begin{array}{l}
T_{p,q,n}^{N}(F_{n})(x_{q},\ldots ,x_{p}) \\
\\
:=\int \mathcal{Q}_{p,n}(x_{p},d(x_{p+1},\ldots ,x_{n}))~\mathcal{M}%
_{q}^{N}(x_{q},d(x_{0},\ldots ,x_{q-1}))~F_{n}(x_{0},\ldots ,x_{n})%
\end{array}%
\end{equation*}

Using formula (\ref{lemformu}), we prove that for any $m\geq 1$ and any
function $F$ on $E_{[q,p]}$
\begin{equation*}
\sqrt{N}~\mathbb{E}\left( \left\vert \left[ \mathcal{M}_{p,q,\eta _{q}^{N}}-%
\mathcal{M}_{p,q,\Phi _{q}\left( \eta _{q-1}^{N}\right) }\right] \left(
F\right) (x_{p})\right\vert ^{m}~\left\vert ~\mathcal{F}_{q-1}^{N}\right.
\right) ^{\frac{1}{m}}\leq a(m)~b(n)~\Vert F\Vert ~\left( \frac{h_{p}^{+}}{%
h_{p}^{-}}(x_{p})\right) ^{2}
\end{equation*}%
for some finite constants $a(m)<\infty $ and $b(n)<\infty $ whose values
only depend on the parameters $m$ and $n$. Using these almost sure
estimates, we easily prove (\ref{keyest}). This ends the proof of the
proposition.\hfill \hbox{\vrule height 5pt width 5pt depth 0pt}\medskip
\newline

Now, we come to the proof of theorem~\ref{tclf}.

\textbf{Proof of theorem~\ref{tclf}:}

Using theorem~\ref{theo2}, we have the decomposition
\begin{equation*}
W_{n}^{\Gamma ,N}(F_{n})=\sum_{p=0}^{n}\gamma _{p}^{N}(1)~V_{p}^{N}\left(
D_{p,n}(F_{n})\right) +R_{n}^{\Gamma ,N}(F_{n})
\end{equation*}%
with the second order remainder term
\begin{equation*}
R_{n}^{\Gamma ,N}(F_{n}):=\sum_{p=0}^{n}\gamma _{p}^{N}(1)~V_{p}^{N}\left(
F_{p,n}^{N}\right) \quad \mbox{\rm and the function}\quad
F_{p,n}^{N}:=[D_{p,n}^{N}-D_{p,n}](F_{n})
\end{equation*}%
By Slutsky's lemma and by the continuous mapping theorem it clearly suffices
to check that $R_{n}^{\Gamma ,N}(F_{n})$ converge to $0$, in probability, as
$N\rightarrow \infty $. To prove this claim, we notice that
\begin{equation*}
\mathbb{E}\left( V_{p}^{N}\left( F_{p,n}^{N}\right) ^{2}~\left\vert ~%
\mathcal{F}_{p-1}^{N}\right. \right) \leq \Phi _{p}\left( \eta
_{p-1}^{N}\right) \left( \left( F_{p,n}^{N}\right) ^{2}\right)
\end{equation*}%
On the other hand, we have
\begin{equation*}
\begin{array}{l}
\Phi _{p}\left( \eta _{p-1}^{N}\right) \left( \left( F_{p,n}^{N}\right)
^{2}\right) \\
\\
=\int \lambda _{p}(dx_{p})~\Psi _{G_{p-1}}\left( \eta _{p-1}^{N}\right)
\left( H_{p}(\mbox{\LARGE .},x_{p})\right) ~F_{p,n}^{N}(x_{p})^{2} \\
\\
\leq \eta _{p}\left( \left( F_{p,n}^{N}\right) ^{2}\right) +\displaystyle%
\int \lambda _{p}(dx_{p})~\left\vert \left[ \Psi _{G_{p-1}}\left( \eta
_{p-1}^{N}\right) -\Psi _{G_{p-1}}\left( \eta _{p-1}\right) \right] \left(
H_{p}(\mbox{\LARGE .},x_{p})\right) \right\vert ~F_{p,n}^{N}(x_{p})^{2}%
\end{array}%
\end{equation*}%
This yields the rather crude estimate
\begin{equation*}
\begin{array}{l}
\Phi _{p}\left( \eta _{p-1}^{N}\right) \left( \left( F_{p,n}^{N}\right)
^{2}\right) \\
\\
=\int \lambda _{p}(dx_{p})~\Psi _{G_{p-1}}\left( \eta _{p-1}^{N}\right)
\left( H_{p}(\mbox{\LARGE .},x_{p})\right) ~F_{p,n}^{N}(x_{p})^{2} \\
\\
\leq \eta _{p}\left( \left( F_{p,n}^{N}\right) ^{2}\right) +4\Vert
Q_{p,n}(1)\Vert ^{2}\displaystyle\int \lambda _{p}(dx_{p})~\left\vert \left[
\Psi _{G_{p-1}}\left( \eta _{p-1}^{N}\right) -\Psi _{G_{p-1}}\left( \eta
_{p-1}\right) \right] \left( H_{p}(\mbox{\LARGE .},x_{p})\right) \right\vert%
\end{array}%
\end{equation*}%
from which we conclude that
\begin{equation*}
\begin{array}{l}
\mathbb{E}\left( V_{p}^{N}\left( F_{p,n}^{N}\right) ^{2}\right) \\
\\
\leq \int \eta _{p}(dx_{p})~\mathbb{E}\left[ \left(
F_{p,n}^{N}(x_{p})\right) ^{2}\right] ~ \\
\\
\qquad +4\Vert Q_{p,n}(1)\Vert ^{2}\int \lambda _{p}(dx_{p})~\mathbb{E}%
\left( \left\vert \left[ \Psi _{G_{p-1}}\left( \eta _{p-1}^{N}\right) -\Psi
_{G_{p-1}}\left( \eta _{p-1}\right) \right] \left( H_{p}(\mbox{\LARGE .}%
,x_{p})\right) \right\vert \right)%
\end{array}%
\end{equation*}%
We can establish that
\begin{equation*}
\sqrt{N}~\mathbb{E}\left( \left\vert \left[ \Psi _{G_{p-1}}\left( \eta
_{p-1}^{N}\right) -\Psi _{G_{p-1}}\left( \eta _{p-1}\right) \right] \left(
H_{p}(\mbox{\LARGE .},x_{p})\right) \right\vert \right) \leq
b(n)~h_{p}^{+}(x_{p})
\end{equation*}%
%
%
%
See for instance section 7.4.3, theorem 7.4.4 in~\cite{fk}. Using
proposition \ref{propDN},
\begin{equation*}
\begin{array}{l}
\sqrt{N}~\mathbb{E}\left( V_{p}^{N}\left( F_{p,n}^{N}\right) ^{2}\right)
\leq c(n)\left( \frac{1}{\sqrt{N}}~\eta _{p}\left( \left( \frac{h_{p}^{+}}{%
h_{p}^{-}}\right) ^{4}\right) +\lambda _{p}(h_{p}^{+})\right)%
\end{array}%
\end{equation*}%
for some finite constant $c(n)<\infty $. The end of the proof of the first
assertion now follows standard computations. To prove the second assertion,
we use the following decomposition
\begin{equation*}
\sqrt{N}~[\mathbb{Q}_{n}^{N}-\mathbb{Q}_{n}](F_{n})=\frac{1}{\overline{%
\gamma }_{n}^{N}(1)}~\overline{W}_{n}^{\Gamma ,N}(F_{n}-\mathbb{Q}%
_{n}(F_{n}))
\end{equation*}%
with the random fields $\overline{W}_{n}^{\Gamma ,N}$ defined in (\ref%
{overlW}). We complete the proof using the fact that $\overline{\gamma }%
_{n}^{N}(1)$ tends to $1$, almost surely, as $N\rightarrow \infty $. This
ends the proof of the theorem.\hfill
\hbox{\vrule height 5pt width 5pt depth
0pt}\medskip \newline

We end this section with some comments on the asymptotic variance associated
to the Gaussian fields $W_{n}^{\mathbb{Q}}$. Using (\ref{sgQ}), we prove
that
\begin{equation*}
\mathbb{Q}_{n}=\Psi _{\overline{D}_{p,n}(1)}(\eta _{p})P_{p,n}
\end{equation*}%
with the pair of integral operators $(\overline{D}_{p,n},P_{p,n})$ from $%
\mathcal{B}(E_{[0,n]})$ into $\mathcal{B}(E_{p})$
\begin{equation*}
\overline{D}_{p,n}(F_{n}):=\frac{D_{p,n}(F_{n})}{\eta _{p}Q_{p,n}(1)}=\frac{%
D_{p,n}(1)}{\eta _{p}Q_{p,n}(1)}~P_{p,n}(F_{n})\quad \mbox{\rm and}\quad
P_{p,n}(F_{n}):=\frac{D_{p,n}(F_{n})}{D_{p,n}(1)}
\end{equation*}%
from which we deduce the following formula
\begin{equation}
\begin{array}{l}
\overline{D}_{p,n}(F_{n}-\mathbb{Q}_{n}(F_{n}))(x_{p}) \\
\\
=\overline{D}_{p,n}(1)(x_{p})~\int \left[
P_{p,n}(F_{n})(x_{p})-P_{p,n}(F_{n})(y_{p})\right] ~\Psi _{\overline{D}%
_{p,n}(1)}(\eta _{p})(dy_{p})%
\end{array}
\label{eqreff}
\end{equation}%
Under condition $(M)_{m}$, for any function $F_{n}$ with oscillations $%
\mbox{\rm
osc}(F_{n})\leq 1$, we prove the following estimate
\begin{equation*}
\Vert \overline{D}_{p,n}(1)\Vert \leq \delta ^{m}\rho \Longrightarrow
\mathbb{E}\left( W_{n}^{\mathbb{Q}}(F_{n})^{2}\right) \leq (\delta ^{m}\rho
)^{2}\sum_{p=0}^{n}\beta (P_{p,n})^{2}
\end{equation*}

\section{Non asymptotic estimates}

\label{nonasymp}

This section is mainly concerned with the proof of theorem~\ref{nonasymptheo}%
. We follow the same semigroup techniques as the ones we used in section
7.4.3 in~\cite{fk} to derive uniform estimates w.r.t. the time parameter for
the $N$-particle measures $\eta _{n}^{N}$. We use the decomposition
\begin{equation*}
\lbrack \mathbb{Q}_{n}^{N}-\mathbb{Q}_{n}](F_{n})=\sum_{0\leq p\leq n}\left(
\frac{\eta _{p}^{N}D_{p,n}^{N}(F_{n})}{\eta _{p}^{N}D_{p,n}^{N}(1)}-\frac{%
\eta _{p-1}^{N}D_{p-1,n}^{N}(F_{n})}{\eta _{p-1}^{N}D_{p-1,n}^{N}(1)}\right)
\end{equation*}%
with the conventions $\eta _{-1}^{N}D_{-1,n}^{N}=\eta _{0}\mathcal{Q}_{0,n}$%
, for $p=0$. Next, we observe that
\begin{equation*}
\begin{array}{l}
\eta _{p-1}^{N}D_{p-1,n}^{N}(F_{n}) \\
\\
=\int \eta _{p-1}^{N}(dx_{p-1})\mathcal{M}_{p-1}^{N}(x_{p-1},d(x_{0},\ldots
,x_{p-2}))\mathcal{Q}_{p-1,n}(x_{p-1},d(x_{p},\ldots
,x_{n}))F_{n}(x_{0},\ldots ,x_{n}) \\
\\
=\int \eta _{p-1}^{N}(dx_{p-1})Q_{p}(x_{p-1},dx_{p}) \\
\\
\hskip3cm\times \mathcal{M}_{p-1}^{N}(x_{p-1},d(x_{0},\ldots ,x_{p-2}))%
\mathcal{Q}_{p,n}(x_{p},d(x_{p+1},\ldots ,x_{n}))F_{n}(x_{0},\ldots ,x_{n})%
\end{array}%
\end{equation*}%
On the other hand, we have
\begin{equation*}
\eta _{p-1}^{N}(dx_{p-1})Q_{p}(x_{p-1},dx_{p})=\eta
_{p-1}^{N}Q_{p}(dx_{p})~M_{p,\eta _{p-1}^{N}}(x_{p},dx_{p-1})
\end{equation*}%
from which we conclude that
\begin{equation*}
\eta _{p-1}^{N}D_{p-1,n}^{N}(F_{n})=(\eta
_{p-1}^{N}Q_{p})(D_{p,n}^{N}(F_{n}))
\end{equation*}%
This yields the decomposition
\begin{equation}
\lbrack \mathbb{Q}_{n}^{N}-\mathbb{Q}_{n}](F_{n})=\sum_{0\leq p\leq n}\left(
\frac{\eta _{p}^{N}D_{p,n}^{N}(F_{n})}{\eta _{p}^{N}D_{p,n}^{N}(1)}-\frac{%
\Phi _{p}(\eta _{p-1}^{N})(D_{p,n}^{N}(F_{n}))}{\Phi _{p}(\eta
_{p-1}^{N})(D_{p,n}^{N}(1))}\right)  \label{fineq}
\end{equation}%
with the convention $\Phi _{0}(\eta _{-1}^{N})=\eta _{0}$, for $p=0$. If we
set
\begin{equation*}
\widetilde{F}_{p,n}^{N}=F_{n}-\frac{\Phi _{p}(\eta
_{p-1}^{N})(D_{p,n}^{N}(F_{n}))}{\Phi _{p}(\eta _{p-1}^{N})(D_{p,n}^{N}(1))}
\end{equation*}%
then every term in the r.h.s. of (\ref{fineq}) takes the following form
\begin{equation*}
\frac{\eta _{p}^{N}D_{p,n}^{N}(\widetilde{F}_{p,n}^{N})}{\eta
_{p}^{N}D_{p,n}^{N}(1)}=\frac{\eta _{p}Q_{p,n}(1)}{\eta _{p}^{N}Q_{p,n}(1)}%
\times \left[ \eta _{p}^{N}\overline{D}_{p,n}^{N}(\widetilde{F}%
_{p,n}^{N})-\Phi _{p}(\eta _{p-1}^{N})\overline{D}_{p,n}^{N}(\widetilde{F}%
_{p,n}^{N})\right]
\end{equation*}%
with the integral operators $\overline{D}_{p,n}^{N}$ defined in (\ref%
{overlDN}). Next, we observe that $D_{p,n}^{N}(1)=Q_{p,n}(1)$, and $%
\overline{D}_{p,n}^{N}(1)=\overline{D}_{p,n}(1)$. Thus, in terms of the
local sampling random fields $V_{p}^{N}$, we have proved that
\begin{equation}
\frac{\eta _{p}^{N}D_{p,n}^{N}(\widetilde{F}_{p,n}^{N})}{\eta
_{p}^{N}D_{p,n}^{N}(1)}=\frac{1}{\sqrt{N}}\times \frac{1}{\eta _{p}^{N}%
\overline{D}_{p,n}(1)}\times V_{p}^{N}\overline{D}_{p,n}^{N}(\widetilde{F}%
_{p,n}^{N})  \label{locaterms}
\end{equation}%
and
\begin{equation}
\overline{D}_{p,n}^{N}(F_{n})=\overline{D}_{p,n}(1)\times
P_{p,n}^{N}(F_{n})\quad \mbox{\rm with}\quad P_{p,n}^{N}(F_{n}):=\frac{%
D_{p,n}^{N}(F_{n})}{D_{p,n}^{N}(1)}\   \label{defPpnN}
\end{equation}%
>From these observations, we prove that
\begin{equation*}
\frac{\Phi _{p}(\eta _{p-1}^{N})(D_{p,n}^{N}(F_{n}))}{\Phi _{p}(\eta
_{p-1}^{N})(D_{p,n}^{N}(1))}=\frac{\Phi _{p}(\eta
_{p-1}^{N})(Q_{p,n}(1)~P_{p,n}^{N}(F_{n}))}{\Phi _{p}(\eta
_{p-1}^{N})(Q_{p,n}(1))}=\Psi _{Q_{p,n}(1)}\left( \Phi _{p}(\eta
_{p-1}^{N})\right) P_{p,n}^{N}(F_n)
\end{equation*}%
Arguing as in (\ref{eqreff}) we obtain the following decomposition
\begin{equation*}
\begin{array}{l}
\overline{D}_{p,n}^{N}(\widetilde{F}_{p,n}^{N})(x_{p}) \\
\\
=\overline{D}_{p,n}(1)(x_{p})\times \int ~\left[
P_{p,n}^{N}(F_{n})(x_{p})-P_{p,n}^{N}(F_{n})(y_{p})\right] ~\Psi
_{Q_{p,n}(1)}(\Phi _{p}(\eta _{p-1}^{N}))(dy_{p})%
\end{array}%
\end{equation*}%
and therefore
\begin{equation*}
\left\Vert \overline{D}_{p,n}^{N}(\widetilde{F}_{p,n}^{N})\right\Vert \leq
b_{p,n}~\beta (P_{p,n}^{N})~\mbox{\rm osc}(F_{n})\quad \mbox{\rm with}\quad
b_{p,n}\leq \sup_{x_{p},y_{p}}{\frac{Q_{p,n}(1)(x_{p})}{Q_{p,n}(1)(y_{p})}}~
\end{equation*}%
We end the proof of (\ref{LLr}) using the fact that for any $r\geq 1$, $%
p\geq 0$, $f\in \mathcal{B}(E_{p})$ s.t. $\mbox{\rm osc}(f)\leq 1$ we have
the almost sure Kintchine type inequality
\begin{equation*}
\mathbb{E}\left( \left\vert V_{p}^{N}(f)\right\vert ^{r}~\left\vert ~%
\mathcal{F}_{p-1}^{N}\right. \right) ^{\frac{1}{r}}\leq a_{r}
\end{equation*}%
for some finite (non random) constants $a_{r}<\infty $ whose values only
depend on $r$. Indeed, using the fact that each term in the sum of (\ref%
{fineq}) takes the form (\ref{locaterms}) we prove that
\begin{equation*}
\sqrt{N}~\mathbb{E}\left( \left\vert [\mathbb{Q}_{n}^{N}-\mathbb{Q}%
_{n}](F_{n})\right\vert ^{r}\right) ^{\frac{1}{r}}\leq a(r)~\sum_{0\leq
p\leq n}~b_{p,n}^{2}~\mathbb{E}\left( \mbox{\rm osc}(P_{p,n}^{N}(F_{n}))%
\right)
\end{equation*}%
This ends the proof of the first assertion (\ref{LLr}) of theorem~\ref%
{nonasymptheo}. For linear functionals of the form (\ref{additivefunct}), it
is easily checked that
\begin{equation*}
D_{p,n}^{N}(F_{n})=Q_{p,n}(1)~\sum_{0\leq q\leq p}\left[ M_{p,\eta
_{p-1}^{N}}\ldots M_{q+1,\eta _{q}^{N}}\right] (f_{q})+\sum_{p<q\leq
n}Q_{p,q}(f_{q}~Q_{q,n}(1))
\end{equation*}%
with the convention $M_{p,\eta _{p-1}^{N}}\ldots M_{p+1,\eta _{p}^{N}}=Id$,
the identity operator, for $q=p$. Recalling that $D_{p,n}^{N}(1)=Q_{p,n}(1)$%
, we conclude that
\begin{equation*}
P_{p,n}^{N}(F_{n})=f_{p}+\sum_{0\leq q<p}\left[ M_{p,\eta _{p-1}^{N}}\ldots
M_{q+1,\eta _{q}^{N}}\right] (f_{q})+\sum_{p<q\leq n}\frac{%
Q_{p,q}(Q_{q,n}(1)~f_{q})}{Q_{p,q}(Q_{q,n}(1))}
\end{equation*}%
and therefore
\begin{equation*}
P_{p,n}^{N}(F_{n})=\sum_{0\leq q<p}\left[ M_{p,\eta _{p-1}^{N}}\ldots
M_{q+1,\eta _{q}^{N}}\right] (f_{q})+\sum_{p\leq q\leq n}\frac{%
Q_{p,q}(Q_{q,n}(1)~f_{q})}{Q_{p,q}(Q_{q,n}(1))}
\end{equation*}%
\begin{equation*}
\frac{Q_{p,q}(Q_{q,n}(1)~f_{q})}{Q_{p,q}(Q_{q,n}(1))}=\frac{S_{p,q}(%
\overline{Q}_{q,n}(1)~f_{q})}{S_{p,q}(\overline{Q}_{q,n}(1))}\quad
\mbox{\rm
with}\quad S_{p,q}(g)=\frac{Q_{p,q}(g)}{Q_{p,q}(1)}
\end{equation*}%
with the potential functions $\overline{Q}_{q,n}(1)$ defined in (\ref%
{defqover}). After some elementary computations, we obtain the following
estimates
\begin{equation*}
\begin{array}{l}
\mbox{\rm osc}(P_{p,n}^{N}(F_{n})) \\
\\
\leq \sum_{0\leq q<p}\beta \left( M_{p,\eta _{p-1}^{N}}\ldots M_{q+1,\eta
_{q}^{N}}\right) \mbox{\rm osc}(f_{q})+\sum_{p\leq q\leq
n}~b_{q,n}^{2}~\beta (S_{p,q})~\mbox{\rm osc}(f_{q})%
\end{array}%
\end{equation*}%
This ends the proof of the second assertion (\ref{estitheo}) of theorem~\ref%
{nonasymptheo}.

\section*{Appendix}

\subsection*{Proof of lemma~\protect\ref{lemMeta}}

We prove the lemma by induction on the parameter $n(>p)$. For $n=p+1$, we
have
\begin{equation*}
\mathcal{M }_{p+1,p,\eta}(x_{p+1},dx_p)=M_{p+1,\eta}(x_{p+1},dx_p)\quad%
\mbox{\rm and}\quad \mathcal{Q }_{p,p+1}(x_p,dx_{p+1})=Q_{p+1}(x_p,dx_{p+1})
\end{equation*}
By definition of the transitions $M_{p+1,\eta}$, we have
\begin{equation*}
\eta Q_{p+1}(dx_{p+1})~ \mathcal{M }_{p+1,p,\eta}(x_{p+1},dx_p)=\eta(dx_p)~
\mathcal{Q }_{p,p+1}(x_p,dx_{p+1})
\end{equation*}
We suppose that the result has been proved at rank $n$. In this situation,
we notice that
\begin{equation*}
\begin{array}{l}
\eta(dx_p)~ \mathcal{Q }_{p,n+1}(x_p,d(x_{p+1},\ldots,x_{n+1})) \\
\\
=\eta(dx_p)~ \mathcal{Q }_{p,n}(x_p,d(x_{p+1},\ldots,x_{n}))
Q_{n+1}(x_n,dx_{n+1}) \\
\\
=\eta Q_{p,n}(dx_n)~Q_{n+1}(x_n,dx_{n+1})~ \mathcal{M }_{n,p,%
\eta}(x_n,d(x_p,\ldots,x_{n-1})) \\
\\
=\eta Q_{p,n}(1)~\Phi_{p,n}(\eta)(dx_n)~Q_{n+1}(x_n,dx_{n+1})~ \mathcal{M }%
_{n,p,\eta}(x_n,d(x_p,\ldots,x_{n-1}))%
\end{array}%
\end{equation*}
Using the fact that
\begin{equation*}
\Phi_{p,n}(\eta)(dx_n)~Q_{n+1}(x_n,dx_{n+1})=
\Phi_{p,n}(\eta)Q_{n+1}(dx_{n+1})~M_{n+1,\Phi_{p,n}(\eta)}(x_{n+1},dx_n)
\end{equation*}
and
\begin{equation*}
\eta Q_{p,n}(1)~\Phi_{p,n}(\eta)Q_{n+1}(dx_{n+1})~=\eta Q_{p,n+1}(dx_{n+1})
\end{equation*}
we conclude that
\begin{equation*}
\begin{array}{l}
\eta(dx_p)~ \mathcal{Q }_{p,n+1}(x_p,d(x_{p+1},\ldots,x_{n+1})) \\
\\
=\eta Q_{p,n+1}(dx_{n+1})~M_{n+1,\Phi_{p,n}(\eta)}(x_{n+1},dx_n) ~ \mathcal{%
M }_{n,p,\eta}(x_n,d(x_p,\ldots,x_{n-1})) \\
\\
=\eta Q_{p,n+1}(dx_{n+1})~ \mathcal{M }_{n+1,p,\eta}\left(x_{n+1},d(x_p,%
\ldots,x_n)\right)%
\end{array}%
\end{equation*}
This ends the proof of the lemma. \hfill%
\hbox{\vrule height 5pt width 5pt
depth 0pt}\medskip \newline

\section*{Proof of lemma~\protect\ref{lemLN}:}

Using the recursions (\ref{recursions}), we prove that
\begin{equation*}
\begin{array}{l}
\mathcal{M}_{n+1,p,\eta _{p}^{N}}\left( x_{n+1},d(x_{p},\ldots ,x_{n})\right)
\\
\\
=\mathcal{M}_{n+1,p+1,\Phi _{p+1}\left( \eta _{p}^{N}\right) }\left(
x_{n+1},d(x_{p+1},\ldots ,x_{n})\right) \times M_{p+1,\eta
_{p}^{N}}(x_{p+1},dx_{p})%
\end{array}%
\end{equation*}%
On the other hand, we also have
\begin{equation*}
\mathcal{M}_{p+1}^{N}(x_{p+1},d(x_{0},\ldots ,x_{p}))=M_{p+1,\eta
_{p}^{N}}(x_{p+1},dx_{p})\mathcal{M}_{p}^{N}(x_{p},d(x_{0},\ldots ,x_{p-1}))
\end{equation*}%
from which we conclude that
\begin{equation*}
\begin{array}{l}
\mathcal{M}_{n+1,p+1,\Phi _{p+1}\left( \eta _{p}^{N}\right) }\left(
x_{n+1},d(x_{p+1},\ldots ,x_{n})\right) \mathcal{M}%
_{p+1}^{N}(x_{p+1},d(x_{0},\ldots ,x_{p})) \\
\\
=\mathcal{M}_{n+1,p,\eta _{p}^{N}}\left( x_{n+1},d(x_{p},\ldots
,x_{n})\right) \mathcal{M}_{p}^{N}(x_{p},d(x_{0},\ldots ,x_{p-1}))%
\end{array}%
\end{equation*}%
The end of the proof is now a direct consequence of the following
decomposition
\begin{equation*}
\begin{array}{l}
\mathcal{M}_{n}^{N}(x_{n},d(x_{0},\ldots ,x_{n-1}))-\mathcal{M}%
_{n}(x_{n},d(x_{0},\ldots ,x_{n-1})) \\
\\
=\sum_{1\leq p\leq n}\left[ \mathcal{M}_{n,p,\eta _{p}^{N}}\left(
x_{n},d(x_{p},\ldots ,x_{n-1})\right) \mathcal{M}_{p}^{N}(x_{p},d(x_{0},%
\ldots ,x_{p-1}))\right. \\
\\
\hskip3cm\left. -\mathcal{M}_{n,p-1,\eta _{p-1}^{N}}\left(
x_{n},d(x_{p-1},\ldots ,x_{n-1})\right) \mathcal{M}%
_{p-1}^{N}(x_{p-1},d(x_{0},\ldots ,x_{p-2}))\right] \\
\\
\hskip.5cm+\mathcal{M}_{n,0,\eta _{0}^{N}}\left( x_{n},d(x_{0},\ldots
,x_{n-1})\right) -\mathcal{M}_{n,0,\eta _{0}}\left( x_{n},d(x_{0},\ldots
,x_{n-1})\right)%
\end{array}%
\end{equation*}%
with the conventions
\begin{equation*}
\mathcal{M}_{n,0,\eta _{0}^{N}}\left( x_{n},d(x_{0},\ldots ,x_{n-1})\right)
\mathcal{M}_{0}^{N}(x_{0},d(x_{0},\ldots ,x_{1}))=\mathcal{M}_{n,0,\eta
_{0}^{N}}\left( x_{n},d(x_{0},\ldots ,x_{n-1})\right)
\end{equation*}%
for $p=0$, and for $p=n$
\begin{equation*}
\mathcal{M}_{n,n,\eta _{n}^{N}}\left( x_{n},d(x_{n},\ldots ,x_{n-1})\right)
\mathcal{M}_{n}^{N}(x_{n},d(x_{0},\ldots ,x_{n-1})=\mathcal{M}%
_{n}^{N}(x_{n},d(x_{0},\ldots ,x_{n-1})
\end{equation*}%
\hfill \hbox{\vrule height 5pt width 5pt depth 0pt}\medskip \newline

\end{document}